%% file: main.tex
\pdfoutput=1
\documentclass[reqno,12pt]{amsart}
\usepackage{a4wide, color, eucal, enumerate, mathrsfs, xcolor}
\usepackage{amsmath,amssymb,epsfig,amsthm}
\usepackage{mathtools}
\usepackage{lipsum}
\usepackage[english]{babel} 
\usepackage[T1]{fontenc}
\usepackage[utf8]{inputenc}
\usepackage[sc]{mathpazo}	
\usepackage{hyperref}
\hypersetup{hidelinks, urlcolor=blue, pdfstartview=FitH, pdfpagemode=UseOutlines,	pdfnewwindow, breaklinks}
\PassOptionsToPackage{unicode}{hyperref}
\usepackage{bookmark}
\usepackage{cleveref}
\usepackage{float}

\input{HeadTikz}

\numberwithin{equation}{section}
\counterwithout{figure}{section}
\newtheorem{theorem}{Theorem}[section]

\newtheorem{proposition}[theorem]{Proposition}

\newtheorem{definition}[theorem]{Definition}

\theoremstyle{remark}
\newtheorem{remark}[theorem]{Remark}

\let\subset\subseteq

\let\phi\varphi
\let\epsilon\varepsilon
\providecommand{\qm}[1]{``#1''}
\newcommand{\wto}{\rightharpoonup}
\DeclareMathOperator{\interior}{int}
\newcommand{\dx}{\, \mathrm{d}x}
\newcommand{\dy}{\, \mathrm{d}y}
\newcommand{\dH}{\, \mathrm{d}\mathcal{H}}
\DeclareMathOperator{\diam}{diam}
\newcommand{\N}{\mathbb{N}}
\newcommand{\R}{\mathbb{R}}
\newcommand{\Z}{\mathbb{Z}}
\def\ue{\mathbf u}
\def\INV{\mathrm{(INV)}}
\DeclareMathOperator{\Deg}{Deg}
\renewcommand{\div}{\operatorname{div}}
\newcommand{\cont}{\mathcal{C}}
\newcommand{\smallMath}[1]{
	\sbox0{$\vcenter{}$}
	\raise\ht0\hbox{%
		\scalebox{0.75}{%
			\lower\ht0\hbox{$#1$}%
		}%
	}
}

\title{Weak limit of $W^{1, 2}$ homeomorphisms in $\R^3$ can have any degree}

\author{Ond\v{r}ej Bouchala}
\address{Czech Technical University in Prague, Faculty of Information Technology, Thákurova 9, 160 00 Prague 6, Czech Republic}
\email{\tt ondrej.bouchala@gmail.com}

\author{Stanislav Hencl}
\address{Department of Mathematical Analysis, Charles University,
So\-ko\-lovsk\'a 83, 186~00 Prague 8, Czech Republic}
\email{\tt hencl@karlin.mff.cuni.cz}

\author{Zheng Zhu}
\address{School of Mathematical Sciences\\
	   Beihang University\\
	  Changping District Shahe Higher Education Park South Third Street No.\ 9\\
	  Beijing 102206\\
	  P.\ R.\ China}
\email{\tt zhzhu@buaa.edu.cn}

\thanks{The first two authors were supported  
by the grant GA\v{C}R P201/21-01976S. The third author was supported by the NSFC grant no.~12301111. This research was done when Z. Zhu was visiting the Department of Mathematical Analysis, Faculty of Mathematics and Physics, Charles University. He wishes to thank Charles University for its hospitality. }
\date{September 2023}

\begin{document}
\begin{abstract}
  {In this paper for every $k\in\Z$ we construct a sequence of weakly converging homeomorphisms $h_m\colon B(0,10)\to\R^3$, $h_m\wto h$ in $W^{1,2}(B(0,10))$, such that $h_m(x)=x$ on $\partial B(0,10)$ and for every $r\in \left(\tfrac5{16},\tfrac{7}{16}\right)$ the degree of $h$ with respect to the ball $B(0,r)$ is equal to $k$ on a set of positive measure.} 
\end{abstract}
\maketitle
\section{Introduction}

In this paper, we study classes of mappings that appear naturally in the modeling of deformations in Continuum Mechanics models. 
Let $\Omega\subset\R^3$ be a domain, i.e.\ a non-empty connected open set, and let $h\colon\Omega\to\R^3$ be a mapping. 
Following the pioneering papers of Ball \cite{Ball} and Ciarlet and Ne\v{c}as \cite{CN} we ask if our mapping is in some sense injective, as the physical \qm{non-interpenetration of the matter} indicates that a deformation should be one-to-one. 

One of the crucial tools in the study of the injectivity of such a mapping (if it is continuous) is the use of topological degree. 
Let $B\subset \Omega$ be a ball, let $h$ be continuous and $y\notin h(\partial B)$. 
Informally speaking the topological degree $\deg(h, B, y)$ is the number of preimages of $y$ under $h$ in $B$ taking the orientation into account (see Preliminaries for more about the degree). For example the topological degree of a homeomorphism is always $0$ or $1$ (if $y\notin h(B)$ or $y\in h(B)$) or always $-1$ or $0$ (if $h$ is reversing the orientation), see e.g. \cite{HM2010} or Proposition \ref{p:top=anal}. 

The topological degree of a continuous mapping depends only on the boundary mapping $h|_{\partial B}$, i.e.\ we can extend it in any continuous way inside and we get the same value of $\deg(h, B,y)$. 
Given $h\in W^{1,p}$, $p>2$, we know that $h|_{\partial B(c, r)}\in W^{1,p}$ for almost all $r$, and hence $h$ has a continuous representative there (as $p>2$ and $2$ is the dimension of the sphere) and we can thus define $\deg(h, B(x,r),y)$. 
This idea was used by M\"uller and Spector \cite{MS1995} (see also e.g.\  \cite{BHMC,HMC, HeMo11,MST,SwaZie2002,SwaZie2004,T}) to define and study the class of mappings in $W^{1,p}(\Omega,\R^3),\ p>2$, that satisfy the so called $\INV$ condition.  
Informally speaking, the $\INV$ condition means that the ball $B(c,r)$ is mapped inside the image of the sphere $h(S(c,r))$ and the complement $\Omega\setminus \overline{B(c,r)}$ is mapped outside $h(S(c,r))$ (see \cite{MS1995} for the formal definition). 
From \cite{MS1995} we know that mappings in this class with $J_h>0$ a.e.\ have   
many desirable properties: they are one-to-one a.e., map disjoint balls into essentially disjoint balls, $\deg(h,B,\cdot)\in\{0,1\}$ for a.e.\ ball $B$, this class is weakly closed and so on.

In all results in the previous paragraph, the authors assume that $h\in W^{1,p}(\Omega,\R^3)$ for some $p>2$. 
However, in some real models for $n=3$, one often works with integrands containing the classical Dirichlet term $|Df|^2$ 
and thus this assumption is too strong. 
Therefore, for $n=3$, Conti and De Lellis \cite{CD2003} introduced the concept of $\INV$ condition also for $W^{1,2}\cap L^{\infty}$ (see also \cite{BHMCR}, \cite{BHMCR2}, \cite{DHM2023} and \cite{DHMo} for some recent work). Note that mappings in $W^{1,2}|_{\partial B(c,r)}$ do not need to be continuous, but we can still define some notion of degree there using the ideas of Brezis and Nirenberg \cite{BN}, i.e.\ we can still define what \qm{inside $h(\partial B(c,r))$} is.  

Mappings satisfying the $\INV$ condition in $W^{1,2}\cap L^{\infty}$ (as defined in \cite{CD2003})  have still all desirable properties but unfortunately this class is not weakly closed and therefore not suitable for the approach of Calculus of Variations. Conti and De Lellis \cite{CD2003} constructed a sequence of bi-Lipschitz mappings $u_m$ in $\R^3$ (which can be easily extended onto $B(0, 10)$ with $u_m(x)=x$ on $\partial B(0,10)$) such that their weak limit $u$ does not satisfy the $\INV$ condition as the degree of the limit is $-1$ on a set of positive measure. This somehow means that the weak limit of $W^{1,2}\cap L^{\infty}$ sense-preserving homeomorphisms can in some sense \qm{revert its orientation} or map something that was inside the ball outside of the image of its sphere. 

The main aim of this paper is to study if the weak limit of $W^{1,2}\cap L^{\infty}$ homeomorphisms can have even more pathological behavior, i.e.\ if the degree could be any $k\in \Z$. Our example shows that this is indeed possible, i.e.\ the weak limit of $W^{1,2}\cap L^{\infty}$ homeomorphisms can somehow \qm{fold itself $k$-times}.

\begin{theorem}\label{main} 
 Let $k\in\Z$. There exists a sequence of $W^{1, 2}$-homeomorphisms 
$\{h_m\}$
from $\overline{B(0, 10)}\subseteq \R^3$ onto $\overline{B(0, 10)}\subseteq \R^3$ such that:
\begin{enumerate}[(i)]
	\item $h_m$ weakly converges to $h$ in $W^{1,2}(B(0,10),\R^3)$.
	\item  $h\big|_{\partial B(0, 10)}$ and $h_m\big|_{\partial B(0, 10)}$ are identity for every $m\in\N$.
\item There exists a set $A\subset B(0, 10)$ with positive measure such that for every $y\in A$ we have 
$$
\Deg \left(h,\ B(0, r),\ y\right)=k
$$
for every $r\in\left(\frac5{16},\frac{7}{16}\right)$. 
\end{enumerate}
\end{theorem}

Such an example of a weak limit of homeomorphisms cannot be constructed in $W^{1,p}$, $p>2$: 
homeomorphisms clearly satisfy the $\INV$ condition and so do their weak limits in $W^{1,p}$, $p>2$, (see \cite[Lemma 3.3]{MS1995}) and mapping in the $\INV$ class with $J_h>0$ a.e. have $\deg(h,B,\cdot)\in\{0,1\}$ for a.e.\ ball $B$ (see \cite[Lemma 3.5]{MS1995}).
Also (under some minor conditions) such an example of a strong limit of homeomorphisms cannot be constructed in  $W^{1,2}$, since the strong limit of $W^{1, 2}$ - homeomorphism also satisfies the $\INV$ condition, please see \cite[Theorem 3.1(b)]{DHM2023}.

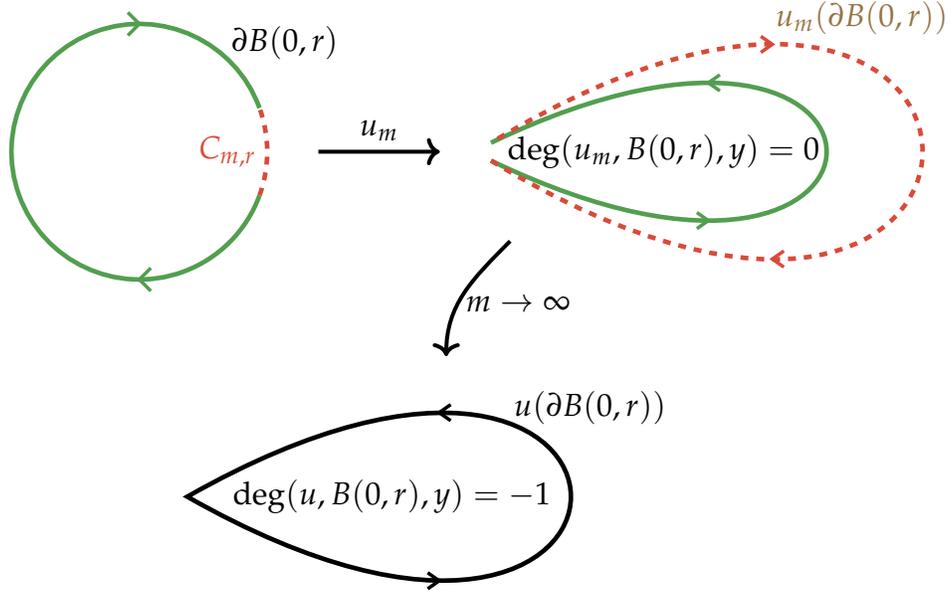
\begin{figure}[H]
	\begin{center}\input{image4.tex}\end{center}
	\caption{Construction of Conti and De Lellis - everything is radially symmetric. Vectors denote the orientation of the sphere and its image. }\label{original}
\end{figure}

To give the reader some idea about the construction of mapping of Conti and De Lellis we show in Fig.\ \ref{original} (everything is radially symmetric around the $x$-axis) how the image of single problematic sphere $\partial B(0,r)$ by $u_m$ and limiting $u$ in the limit looks like. 

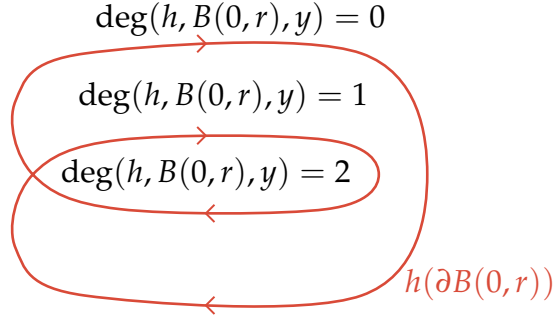
\begin{figure}[H]
	\begin{center}\input{image5.tex}\end{center}
	\caption{Our construction of limit mappings - everything is radially symmetric.}	
	\label{fig:limit}	
\end{figure}

The image of a spherical cap $C_{m,r}\subset \partial B(0,r)$ (of size roughly $1/m$) by $u_m$ is the outer \qm{dashed} bubble (of size roughly $1$) and this bubble disappears in the limit as $\diam(C_{m,r})\overset{m\to\infty}{\to}0$. The limiting mapping has a degree equal to $-1$ inside the inner bubble as the orientation of the bubble is reversed with respect to the original ball. The derivative of $h_m$ on $C_{m,r}$ is comparable to 
$$
\int_{C_{m,r}}|Dh_m|^2\approx \mathcal{H}^2(C_{m,r})\left|\frac{1}{\tfrac{1}{m}}\right|^2\approx \frac{1}{m^2}|m|^2=1
$$ 
and thus it remains uniformly bounded in $m$ even when we integrate it over $r$ and thus $f_m$ have weak limit in $W^{1,2}$. 
Of course, one needs to extend it as a homeomorphism on all neighboring spheres, make it identity on $\partial B(0,10)$, estimate derivatives in all directions but the essential idea is in the figure. 
To obtain degree $k=2$ in our example, we need to achieve something like the following behavior for the limit mapping $h$, see Fig.\ \ref{fig:limit}  (again everything is rotated around the $x$-axis), which is more complicated. We need more than one spherical cap with strange behavior and our $h_m$ are not radially symmetric.

Let us also note that the class of weak limits of Sobolev homeomorphisms was recently characterized in the planar case by Iwaniec and Onninen \cite{IO, IO2} and De Philippis and Pratelli \cite{DPP}. For some other kind of pathological examples of limits of $W^{1,2}\cap L^{\infty}$ homeomorphisms where the limit fails to be one to one a.e., we refer the reader to Bouchala, Hencl and Molchanova \cite{BHM}. For other interesting examples of non-injective mappings see \cite{KKM2001}.


\section{Preliminaries} 
We use $B(x, r)$ for a ball in $\R^3$ with the center $x\in\R^3$ and with the radius $r>0$. By $\interior(A)$ we denote the interior of the set $A$. Throughout the paper, $C$ will be a generic positive constant, which may even be different in a single string of estimates.

\subsection{Degree for continuous mappings}
Let $\Omega\subset\R^3$ be a bounded open set.
Given a continuous map $f\colon\overline\Omega\to\R^3$ and $y\in \R^3\setminus f(\partial\Omega)$,
we can define the {\it topological degree} as
$$\deg(f,\ \Omega,\ y)=\sum_{x\in\Omega\cap f^{-1}(y)} \operatorname{sgn}(J_f(x))$$
if $f$ is smooth in $\Omega$ and $J_f(x)\neq 0$ for each $x\in \Omega\cap f^{-1}(y)$.
By uniform approximation, 
this definition can be extended to an arbitrary continuous mapping 
$f\colon\overline\Omega\to\R^3$. Note that the degree depends only on
values of $f$ on $\partial \Omega$.

If $f\colon \overline\Omega\to\R^n$ is a homeomorphism,
then either $\deg (f,\Omega,y)=1$ for all $y\in f(\Omega)$
($f$ is \textit{sense preserving}), or 
$\deg (f,\Omega,y)=-1$ for all $y\in f(\Omega)$
($f$ is \textit{sense reversing}). If, in addition,
$f\in W^{1, n-1}(\Omega,\R^n)$, then this topological orientation
corresponds to the sign of the Jacobian. More precisely, we have

\begin{proposition}[\cite{HM2010}]\label{p:top=anal} 
  Let $f\in W^{1, n-1}(\Omega,\R^n)$ be 
  a homeomorphism on $\overline\Omega$ with $J_f>0$ a.e.
  Then 
  $$
	  \deg(f,\ \Omega,\ y)=1,\qquad y\in f(\Omega).
  $$
\end{proposition}

\medskip

\subsection{Degree for \texorpdfstring{$W^{1,2}\cap L^{\infty}$}{W\^{}\{1,2\} bounded} mappings}\label{degree}

Let $B$ be a ball,
$f\in W^{1,2}(\partial B,\R^3)\cap \cont(\partial B,\R^3)$, $|f(\partial B)|=0$, 
and $\ue\in \cont^1(\R^3,\R^3)$, then (see \cite[Proposition 2.1]{MS1995})
\begin{equation}\label{weakdegree}
  \int_{\R^3}\deg(f,B,y)\operatorname{div} \ue(y)\dy=
  \int_{\partial B} (\ue\circ f)\cdot (\Lambda_{2} D_{\tau}f)\nu \dH^2,
\end{equation}
where $D_{\tau}f$ denotes the tangential gradient and $\Lambda_{2} D_{\tau}f$ is the restriction of $\operatorname{cof} Df$ to the corresponding subspace (see \cite{DHM2023} for details).

Following \cite{CD2003} (see also \cite{BN}) we need a more general version 
of the degree 
which works for mappings in $W^{1,2}\cap L^{\infty}$ 
that are not necessarily continuous. 

\begin{definition}\label{defdegree}
  Let $B\subset\R^3$ be a ball and let $f\in W^{1,2}(\partial B,\R^3)\cap L^{\infty}(\partial B,\R^3)$. Then we define 
  $\Deg(f, B, \cdot)$ as the distribution satisfying
  \begin{equation}\label{qqq}
	  \int_{\R^3}\Deg(f,B,y)\psi(y)\dy=
	  \int_{\partial B} (\ue\circ f)\cdot(\Lambda_{2} D_{\tau}f) \nu \dH^2
  \end{equation}
  for every test function $\psi\in \cont_c^{\infty}(\R^3)$ 
  and every 
  $\cont^{\infty}$ vector field $\ue$ on $\R^3$ satisfying $\div \ue=\psi$.
\end{definition}

As in \cite{CD2003} (see also \cite{DHM2023}) it can be verified that the right-hand side does not depend on the way 
$\psi$ is expressed as $\div \ue$ and that the distribution $\Deg(f, B,\cdot)$ can be represented as a $BV$ function.

\begin{remark}\label{Deg=deg}
  Let $B$ be a ball and $f\in W^{1,2}(\partial B,\R^3)\cap \cont(\overline B,\R^3)$.
  If $|f(\partial B)|=0$, then $\Deg(f,B,y)=\deg(f,B,y)$ for 
  a.e.\  $y\in\R^3$. 
\end{remark}

\subsection{Matrix of derivatives in different coordinates} 

Let $(x_1,r,\phi)$ denote the usual cylindrical coordinates in $\R^3$ and let $a\colon\R^3\to\R^3$ be a mapping from cylindrical to cylindrical coordinates, i.e.
$$
a(x_1,r,\phi)=(a^{x_1}(x_1,r,\phi),a^r(x_1,r,\phi),a^{\phi}(x_1,r,\phi))
$$
It is well-known that the matrix of derivatives of $a$ in this coordinate system is
\begin{equation}\label{dercyl}
Da(x_1, r, \phi)=\begin{pmatrix}
	\frac{\partial a^{x_1}}{\partial x_1} &
	\frac{\partial a^{x_1}}{\partial r} & 
	\frac1r{\cdot} \frac{\partial a^{x_1}}{\partial\phi} \\[2mm]
	\frac{\partial a^r}{\partial x_1} & 
	\frac{\partial a^r}{\partial r} &
	\frac1r {\cdot}  \frac{\partial a^r}{\partial\phi} \\[2mm]
	\smallMath{a^r} {\cdot}  \frac{\partial a^{\phi}}{\partial x_1} & 
	\smallMath{a^r} {\cdot}   \frac{\partial a^{\phi}}{\partial r} &
	\frac{a^r}{r}   {\cdot} \frac{\partial a^{\phi}}{\partial\phi}
	\end{pmatrix}.
\end{equation}

Let $(r,\theta,\phi)$ denote the usual spherical coordinates in $\R^3$, i.e.\ $\phi\in(0,2\pi)$ and ${\theta\in(0,\pi)}$. Let $b\colon \R^3\to\R^3$ be a mapping from cylindrical to spherical coordinates, i.e.
$$
b(x_1,r,\phi)=(b^r(x_1,r,\phi),b^{\theta}(x_1,r,\phi),b^{\phi}(x_1,r,\phi)). 
$$
It is well-known that the matrix of derivatives of $a$ in this coordinate system is
\begin{equation}\label{dersphr}
Db(x_1, r, \phi)=\begin{pmatrix}
\frac{\partial b^r}{\partial x_1} &
\frac{\partial b^r}{\partial r} &
\frac1r{\cdot}\frac{\partial b^r}{\partial\phi} \\[2mm]
\smallMath{b^r} {\cdot}  \frac{\partial b^\theta}{\partial x_1} &
\smallMath{b^r} {\cdot}  \frac{\partial b^\theta}{\partial r} &
\frac{b^r}{r}{\cdot}\frac{\partial b^\theta}{\partial\phi} \\[2mm]
\smallMath{b^r \sin b^\theta } {\cdot} \frac{\partial b^{\phi}}{\partial x_1} &
\smallMath{b^r \sin b^\theta } {\cdot} \frac{\partial b^{\phi}}{\partial r} &
\frac{b^r\sin b^\theta}{r}{\cdot}\frac{\partial b^{\phi}}{\partial\phi}
\end{pmatrix}. 
\end{equation}

\section{Proof of main theorem for \texorpdfstring{$k=2$}{k=2}}
\subsection{Overview}\ \\
In this subsection we roughly explain the construction of the mappings $h$ and $h_m$ from the Theorem \ref{main}. To create degree two somewhere we need to go around that area twice (imagine the planar case). To achieve that for the limit function $h$, we will define $h_m$ to go around three times, twice in the positive and once in the negative direction (or orientation). Using the idea and the mapping $u_m$ from \cite[Section~6]{CD2003} we can create such loops (or bubbles). And by preparing the set before applying the mapping $u_m$ we can control which bubbles will disappear in the limit. 

\begin{figure}[H]
	\begin{center}\input{image7.tex} \end{center}
	\caption{Overview of the construction.}
	\label{fig:overview}
\end{figure}
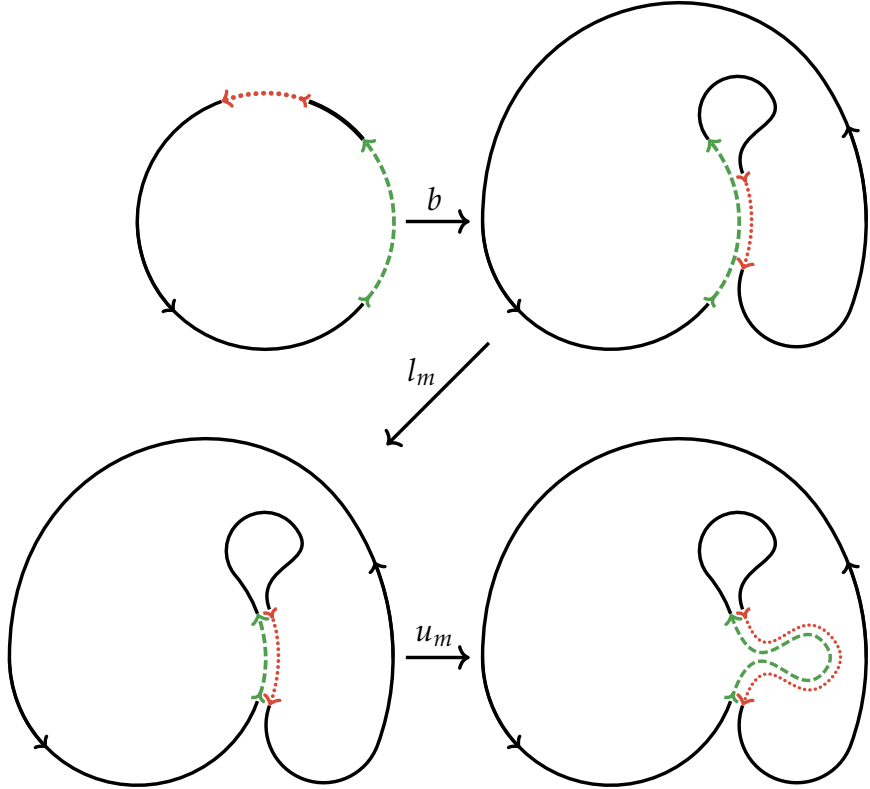

The mappings $h_m$ will be the composition of three homeomorphisms, $$h_m(x):=u_m\circ l_m\circ b(x).$$ The Figure \ref{fig:overview} roughly explains what these three homeomorphisms do. 
The mapping $b$ is bi-Lipschitz and it does not depend on $m$. It moves the dashed part of the circle/sphere (that will not disappear) and the dotted part (which will disappear) close together so that the image of the dashed bubble will be inside the image of the dotted one at the end. 

The mappings $l_m$ are Lipschitz, and they squash the dashed part (of size comparable to $1$) to something smaller (of size comparable to $\tfrac{1}{m}$). In this way, the dashed arc does not disappear in the limit (as it was big at the start) even after applying the Conti -- De Lellis mapping $u_m$. After this careful preparation, we can take the mapping $u_m$ exactly as in the paper \cite{CD2003}. Of course, it would be possible to use also a mapping from \cite[Theorem 1.2]{DHM2023}. 

\subsection{Definition of \texorpdfstring{$b$}{b}}
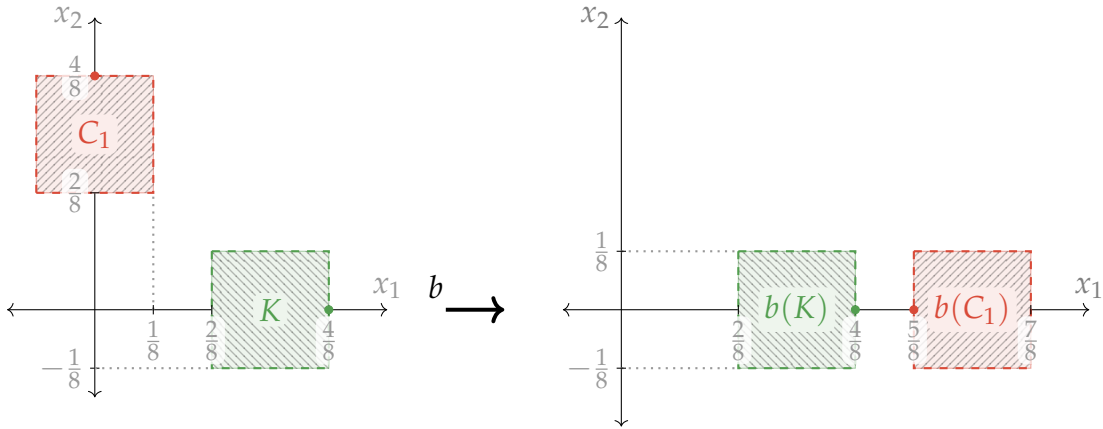
\begin{figure}[H]
 	\begin{center}\input{./image2}\end{center}
	\caption{The mapping $b$ in the plane $\{x_1, x_2\}$.}
	\label{fig:b}
\end{figure}

\ \\
The mapping $b$ is an identity on the set $K$, where 
$$
K:=\left\{(x_1, x_2, x_3)\in B(0, 10):\frac28<x_1<\frac48,\  0\leq\sqrt{x_2^2+x_3^2}<\frac{1}{8}\right\}.\\
$$
We define for every $m\in\mathbb N$ (see Fig.~\ref{fig:b}) 
$$C_m:=\left\{(x_1, x_2, x_3)\in B(0, 10):\frac{2}{8}<x_2<\frac{4}{8},\ 0\leq\sqrt{x_1^2+x_3^2}<\frac1{8m}\right\}.$$
On cylinder ${C_1}$ the mapping $b$ is defined by
\begin{equation*}
  b(x_1, x_2, x_3)=\left(-x_2+\tfrac98, x_1, x_3\right).
  \end{equation*}
That means that on $C_1$ the mapping $b$ is just translation and rotation, see Figure \ref{fig:b}. For every $m\in\N$ it maps $C_m$  onto the set 
$$
b({C_m})=\left\{(x_1, x_2, x_3)\in B(0, 10):\frac58<x_1<\frac78,\  0\leq\sqrt{x_2^2+x_3^2}<\frac1{8m}\right\}.
$$
The Jacobian matrix of $b$ for $x\in C_1$ is 
\begin{equation*}
  Db(x)=
\begin{pmatrix}
0 & -1 & 0 \\
1 & 0 & 0 \\
0 & 0 & 1
\end{pmatrix}.
\end{equation*}
It is not difficult to see that we can extend $b$ to be a bi-Lipschitz homeomorphism which maps $\overline{B(0, 10)}$ onto $\overline{B(0, 10)}$ and which is the identity on the boundary, for example by defining it by hand. We do not give the exact formula here because it would be too lengthy and confusing. 

Later we add additional requirements about the behavior on spheres near $\partial B(0,\tfrac38)$ so that their images do not cross some "forbidden" regions but this can be clearly satisfied as well.

\subsection{Definition of \texorpdfstring{$l_m$}{l\_m}} 
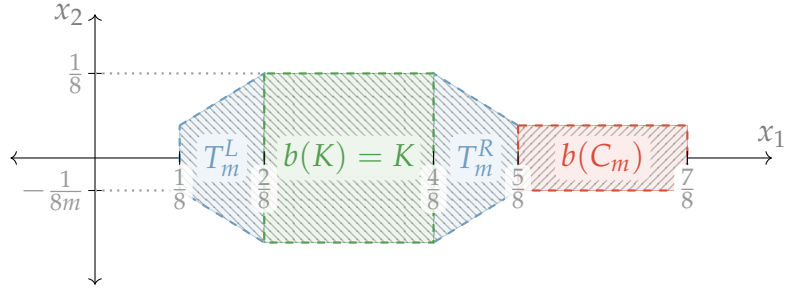
\begin{figure}[H]
	\begin{center}\input{./image1.tex}\end{center}
	\caption{The set $B_m$ in the plane $\{x_1,x_2\}$.}
	\label{fig:Ds}
\end{figure}

\ \\
First, we define a thick cylinder (containing the sets $K=b(K)$ and $b(C_m)$ from earlier) by setting
$$ 
D_0:=\left\{(x_1, x_2, x_3)\in B(0, 10): \tfrac18<x_1<\tfrac78,\  0\leq\sqrt{x_2^2+x_3^2}<\tfrac28\right\}.
$$
The mapping $l_m$ is the identity outside of this cylinder. For every $m\in\mathbb N$, we define a slim cylinder by setting
$$D_m:=\left\{(x_1, x_2, x_3)\in B(0, 10):\tfrac18<x_1<\tfrac78,\  0\leq\sqrt{x_2^2+x_3^2}<\frac{1}{8m}\right\},$$
and we define two cones $T_m^L$ and $T_m^R$ by setting (see Fig.\ \ref{fig:Ds})
\begin{align*}
	T_m^L&:=\left\{(x_1, x_2, x_3)\in B(0, 10): \tfrac18\leq x_1\leq\tfrac28,\ 0\leq\sqrt{x_2^2+x_3^2}<\left(1-\tfrac{1}{m}\right)x_1+\tfrac{2-m}{8m}\right\},\\
T_m^R&:=\left\{(x_1, x_2, x_3)\in B(0, 10):\tfrac48\leq x_1<\tfrac58,\  0\leq\sqrt{x_2^2+x_3^2}<\left(\tfrac{1}{m}-1\right)x_1-\tfrac{5m-4}{8m}\right\}.
\end{align*}
Then, we define a \qm{bonbon-like} domain $B_m$ by setting 
$$B_m:= T_m^L\cup b(K)\cup T_m^R\cup b(C_m) .$$
The mapping $l_m$ squeezes $B_m$ into the slim cylinder $D_m$ and stretches the set $D_0\setminus B_m$ onto  $D_0\setminus D_m$. 

To be precise, we define the homeomorphism $l_m$ inside $D_0$ using the cylindrical coordinates $(x_1, r, \phi)$, that is $(x_1,x_2,x_3)=(x_1, r\cos(\phi), r\sin(\phi))$:

The definition of $l_m$ for $x=(x_1,r,\phi)\in B_m$ is
$$
l_m(x) := \begin{cases}
	\left(x_1,\ \frac{r}{(8m-8)x_1+2-m},\ \phi\right),&x\in T_m^L,\\[0.2cm]
	\left(x_1,\ \tfrac{r}{m},\ \phi\right),&x\in b(K),\\[0.2cm]
	\left(x_1,\ \frac{r}{(8-8m)x_1+5m-4},\ \phi\right),&x\in T_m^R,\\[0.2cm]
	x,&x\in b(C_m).
\end{cases}
$$
On the set $D_0\setminus{B_m}$ we define $l_m$ to be linear, with respect to the radius $r$ (and mapping $x_1$ to $x_1$ and $\phi$ to $\phi$), such that it maps the set $D_0\setminus B_m$ onto the \qm{annulus} $D_0\setminus D_m$ . It is not difficult to see that $l_m$ is Lipschitz on $D_0\setminus B_m$. 

Next, we compute the matrix of derivatives of $l_m$. Obviously, this matrix is the identity on $C_m$. The matrix of derivatives of $l_m$ on $T_m^L$, with respect to the cylindrical coordinates $(x_1, r, \phi)$, is (see \eqref{dercyl})
\begin{equation}\label{eq:matrix1}
	Dl_m(x)=
	\begin{pmatrix}
	1 & 0 & 0\\
	\frac{(8-8m)r}{\left((8m-8)x_1+2-m\right)^2} & \frac{1}{(8m-8)x_1+2-m} & 0\\
	0 & 0 & \tfrac{1}{(8m-8)x_1+2-m}
	\end{pmatrix}.
	\end{equation}

The matrix of derivatives of $l_m$ on $b(K)$ with respect to the cylindrical coordinates is
\begin{equation}\label{eq:matrix2}
Dl_m(x)=
\begin{pmatrix}
1 & 0 & 0\\
0 & \tfrac{1}{m} & 0\\
0 & 0 & \tfrac{1}{m}
\end{pmatrix}.
\end{equation}

And on $T_m^R$ the matrix of derivatives w.r.t. the cylindrical coordinates is 
\begin{equation}\label{eq:matrix3}
Dl_m(x)=
\begin{pmatrix}
1 & 0 & 0\\
\frac{(8m-8)r}{\left((8-8m)x_1+5m-4\right)^2} & \frac{1}{(8-8m)x_1+5m-4} & 0\\
0 & 0 & \tfrac{1}{(8-8m)x_1+5m-4}
\end{pmatrix}.
\end{equation}

\subsection{Definition of \texorpdfstring{$u_m$}{u\_m}}
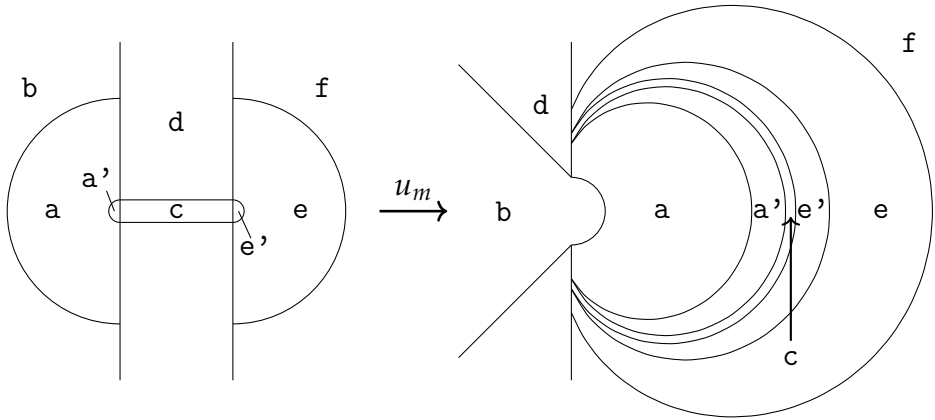
\begin{figure}[H]
	\begin{center}\input{image6.tex}\end{center}	
	\caption{The areas as in \cite[Figure 2.]{CD2003}.}\label{obrCD}
\end{figure}

\ \\
The mapping $u_m$ is the same as in \cite[Section 6]{CD2003} with $\epsilon=\tfrac1{8m}$. They define it piecewise in several regions. For convenience, we include their picture here as Fig.\ \ref{obrCD}. The mapping is axially symmetric with respect to $x_1$-axis so Fig.\ \ref{obrCD} is rotated around it and the regions \texttt{a, a', b, c, d, e, e'} on the left are mapped by the mapping to the corresponding regions on the right. In the limit the regions \texttt{a', c, e'} disappear and the region \texttt{a} which was outside of \texttt{e} is mapped inside the image of the boundary of \texttt{e}. 

For us it is important how it works in the region that contains our bonbons $B_m$, that is in the region \texttt{c},
$$ \texttt{c}:=\left\{(x_1, r, \phi)\in \R^3: 0\leq x_1\leq 1,\ 0\leq r\leq\tfrac{1}{8m}\right\}. $$
The mapping $u_m$ in the region \texttt{c} is defined from cylindrical to spherical coordinates as
\begin{equation}
  \label{eq:vm}
  \begin{aligned}
	  u_m(x_1, r, \phi)&:=\left(u_m^r(x_1,r, \phi), u_m^\theta(x_1,r, \phi), \phi\right),\\
	  u_m^r(x_1,r, \phi)&:=\left(1+\frac{1}{8m}\right)\cos\left(2\arctan (8m\cdot r)\right)+\frac{2}{8m}+\frac{x_1}{64m^2},\\
	  u_m^\theta(x_1,r, \phi)&:=2\arctan(8m\cdot r).
  \end{aligned}
\end{equation}

Obviously $l_m(B_m)$ is a subset of \texttt{c}. Since $u_m^r$ does not depend on $\phi$ and $u_m^\theta$ does not depend on $x_1$ or $\phi$, the matrix of derivatives of $u_m$ at the point $(x_1, r, \phi)\in\texttt{c}$ with respect to the cylindrical coordinates is (see \eqref{dersphr})
\begin{equation}\label{eq:Javm}
Du_m(x_1,r, \phi)=
\begin{pmatrix}
\frac{\partial u_m^r}{\partial x_1} & \frac{\partial u_m^r}{\partial r} & 0\\
0 & u_m^r{\cdot}\frac{\partial u_m^\theta}{\partial r} & 0\\
0 & 0 & \frac{u_m^r\sin(u_m^\theta)}{r} 
\end{pmatrix}.
\end{equation}
It is easy to check that 
\begin{align*}
  \frac{\partial u_m^r}{\partial x_1}&=\frac{1}{64m^2},\\
  \frac{\partial u_m^r}{\partial r}&=-\left(1+\frac{1}{8m}\right)\sin\left(2\arctan(8m\cdot r)\right)\frac{16m}{1+64m^2r^2},\\
  u_m^r\frac{\partial u_m^\theta}{\partial r}&=
  \left(\left(1+\frac{1}{8m}\right)\cos(2\arctan(8m\cdot r))+\frac{1}{4m}+\frac{x_1}{64m^2}\right)\frac{16m}{1+64m^2r^2},\\
  \frac{u_m^r\sin(u_m^\theta)}{r}&=\left(\left(1+\frac{1}{8m}\right)\cos(2\arctan(8m\cdot r))+\frac{1}{4m}+\frac{x_1}{64m^2}\right)\frac{\sin(2\arctan(8m\cdot r))}{r}.
\end{align*}

For $x=(x_1, r, \theta)\in\texttt{c}$ there is a positive constant $C$ independent on $m$ such that
\begin{equation}\label{eq:bound}
\begin{array}{ccc}
	\phantom{m\cdot}0\leq \left|\dfrac{\partial u_m^r(x)}{\partial x_1}\right|\leq \frac{C}{m^2}\phantom{\cdot m}&&
	\phantom{m\cdot}0\leq\left|\dfrac{\partial u_m^r(x)}{\partial r}\right|\leq C \cdot m,\\[5mm]
	  \phantom{m\cdot}0\leq\left|u_m^r(x)\dfrac{\partial u_m^\theta(x)}{\partial r}\right|\leq C\cdot   m,&&
	\phantom{m\cdot}0\leq\left|\dfrac{u_m^r(x)\sin(u_m^\theta(x))}{r}\right|\leq C\cdot   m.
\end{array}
\end{equation}
  
In their paper \cite{CD2003} they do not need the mappings $u_m$ to be identity on the boundary. But it is not difficult to observe that their mappings are bi-Lipschitz and \qm{well-behaved} on $\partial B(0,3)$, so we can extend them by hand to be uniformly bi-Lipschitz on $B(0,10)\setminus B(0, 3)$ and identity on $\partial B(0,10)$. So we get a sequence of homeomorphisms with

\begin{equation}\label{eq:1}
\sup_{m\in\N}\int_{B(0, 10)}\left|D u_m(y)\right|^2\dy<\infty.
\end{equation}

\begin{figure}[H]
	\begin{center}  \input{./image3}  \end{center}
 	\caption{The mapping $b$ near the sphere $\partial B(0,\tfrac38)$ in the plane $\{x_1,x_2\}$.}
	\label{fig:wigglycircle}
\end{figure}
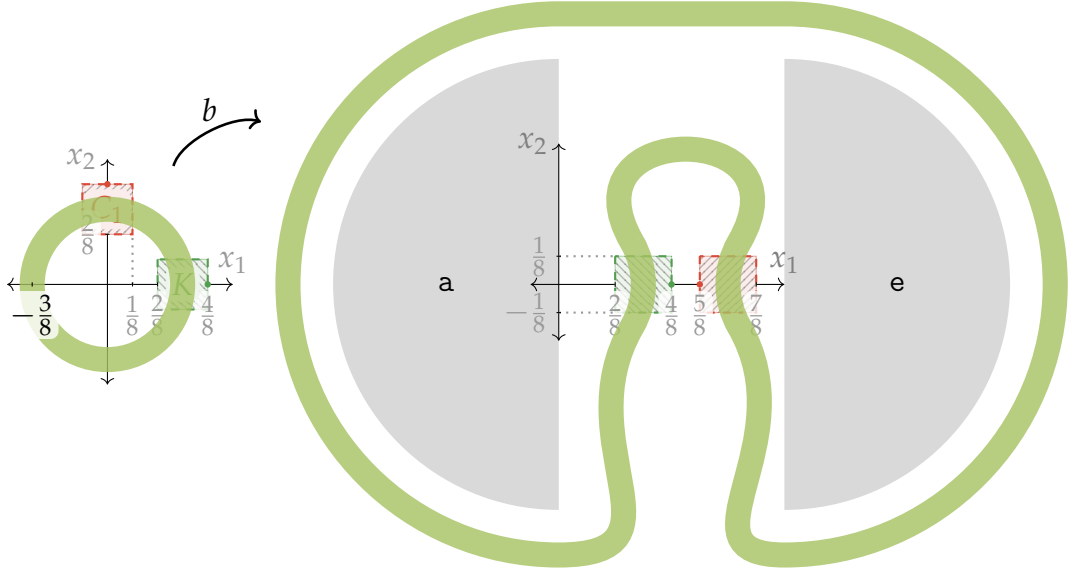

As we have already mentioned we need a small additional requirement about the behavior of our first map $b$ on the set $A:=B(0,\tfrac{7}{16})\setminus B(0,\tfrac{5}{16})$. 
We need $A$ to be mapped to the regions \texttt{b}, \texttt{d} and \texttt{f} from the Figure \ref{obrCD} on the left (i.e.\ they do not enter $a$, $a'$, $e$ and $e'$), except for the intersection of $A$ with $b(K)$ and $b(C_1)$, which can be mapped into the set \texttt{c} from the same figure. This could be easily achieved as seen in Fig.\ \ref{fig:wigglycircle}. 

We shall also need that $l_m(D_0\setminus B_m)$ is the subset of \texttt{d} (see Fig.\ \ref{fig:wigglycircle}) and that the derivative of $u_m$ is bounded by a constant (independent on $m$) as $u_m$ are uniformly Lipchitz there (see \cite[construction of \texttt{d} - pages 544-545]{CD2003}).  

\subsection{Computation of the \texorpdfstring{$\mathbf{W^{1,2}}$}{W\^{}\{1,2\}} norm.}\ \\
Now, we set
$$
h_m(x):=u_m\circ l_m\circ b(x).
$$
Since all three mappings are homeomorphisms it is easy to see that $h_m$ is homeomorphism as well. Analogously it is easy to see that $h_m(x)=x$ for $x\in\partial B(0,10)$.

Since  $b$  is bi-Lipschitz, in order to show 
$$\sup_{m\in\N}\int_{B(0, 10)}\left|Dh_m(x)\right|^2\dx<\infty,$$
it suffices to show
$$\sup_{m\in\N}\int_{B(0, 10)}\left|D(u_m\circ l_m)(x)\right|^2\dx<\infty.$$
Since  $l_m$ is identity on $B(0, 10)\setminus D_0$, from (\ref{eq:1}), we obtain 
\begin{equation}\label{eq:2}
\sup_{m\in\N}\int_{B(0, 10)\setminus D_0}\left|D(u_m\circ l_m)(x)\right|^2\dx<\infty.
\end{equation}
As mentioned at the end of the previous subsection $u_m$ are uniformly Lipschitz on $l_m( D_0\setminus B_m)$. Hence, for every $x\in D_0\setminus\overline B_m$, we have 
$$|D(u_m\circ l_m)(x)|\leq |Du_m(l_m(x))|\cdot|Dl_m(x)|\leq C.$$
Therefore
\begin{equation}\label{eq:3}
   \sup_{m\in\N}\int_{ D_0\setminus\overline B_m}\left|D(u_m\circ l_m)(x)\right|^2\dx<\infty.
\end{equation}

It remains to consider the derivative on $B_m$. 
For almost every $x\in B_m$ the chain rule (with respect to the correct system of coordinates) gives 
\begin{equation}\label{eq:trans}
D(u_m\circ l_m)(x)=Du_m(l_m(x))\cdot Dl_m(x).
\end{equation}
Since $l_m$ is always identity on $C_m$, by (\ref{eq:1}), we have
\begin{equation}\label{eq:4}
\sup_{m}\int_{C_m}|D(u_m\circ l_m)(x)|^2\dx<\infty.
\end{equation}
By (\ref{eq:matrix2}), (\ref{eq:Javm}) and (\ref{eq:trans}) we know that for every $x\in b(K)$
$$ D(u_m\circ l_m)(x)=
\begin{pmatrix}
\frac{\partial u_m^r}{\partial x_1} & \frac1m{\cdot}\frac{\partial u_m^r}{\partial r} & 0\\
0 & \frac{u_m^r}{m}{\cdot}\frac{\partial u_m^\theta}{\partial r} & 0\\
0 & 0 & \frac{u_m^r\sin(u_m^\theta)}{mr}
\end{pmatrix}$$
and (\ref{eq:bound}) now imply $\left|D(u_m\circ l_m)(x)\right|\leq C .$ Hence we have 
\begin{equation}\label{eq:5}
\sup_{m\in\N}\int_{b(K)}\left|D(u_m\circ l_m)(x)\right|^2\dx<\infty.
\end{equation}

Let 
$$\diamondsuit:=(8m-8)x_1+2-m.$$
By (\ref{eq:matrix1}), (\ref{eq:Javm}) and (\ref{eq:trans}), for every $x\in\interior(T_m^L)$ it holds that 
\begin{equation}\label{matrix}
D(u_m\circ l_m)(x)=
\begin{pmatrix}
	\frac{\partial u_m^r}{\partial x_1}+\frac{(8-8m)r}{\diamondsuit^2}{\cdot}\frac{\partial u_m^r}{\partial r}  &
	 \tfrac{1}{\diamondsuit}{\cdot}\frac{\partial u_m^r}{\partial r}&
	  0\\[2mm]
\frac{(8-8m)r \cdot u_m^r}{\diamondsuit^2}{\cdot}\frac{\partial u_m^\theta}{\partial r} &
 \frac{u_m^r}{\diamondsuit}{\cdot}\frac{\partial u_m^\theta}{\partial r} &
  0\\[2mm]
0 & 0 & \frac{u_m^r\sin(u_m^\theta)}{r\cdot \diamondsuit}
\end{pmatrix}.
\end{equation}
By definition of $T_m^L$ we know that for $x\in \interior(T_m^L)$ we have 
\begin{equation}\label{eq:diamond}
0\leq r<\left(1-\frac{1}{m}\right)x_1+\frac{2-m}{8m}=\frac{\diamondsuit}{8m},
\end{equation}
so 
$$\left|\frac{(8-8m)r}{\diamondsuit}\right|\leq C$$
for a constant $C$ independent on $x_1$, $r$ and $m$. By \eqref{matrix} and (\ref{eq:bound}), for every $x\in\interior(T_m^L)$, we have 
$$\left|D(u_m\circ l_m)(x)\right|\leq\frac{C\cdot m}{\diamondsuit}.$$
Hence the Fubini theorem implies (see also \eqref{eq:diamond})
\begin{equation}\label{eq:6}
\sup_{m\in\N}\int_{\interior(T_m^L)}\left|D(u_m\circ l_m)(x)\right|^2\dx 
\leq C\int_{\frac{1}{8}}^{\frac{2}{8}}\left(\frac{\diamondsuit}{8m}\right)^2\frac{m^2}{\diamondsuit^2}\dx_1<\infty.
\end{equation}
Since $T_m^R$ is essentially the same as $T_m^L$, similar computation gives
\begin{equation}\label{eq:7}
\sup_{m\in\N}\int_{\interior( T_m^R)}\left|D(u_m\circ l_m)(x)\right|^2\dx<\infty.
\end{equation}
The boundaries of $ D_0$, $C_m$, $b(K)$, $T_m^L$ and $T_m^R$ have zero measure, so after summing the inequalities (\ref{eq:2}), (\ref{eq:3}), (\ref{eq:4}), (\ref{eq:5}), (\ref{eq:6}) and (\ref{eq:7}) we obtain the desired inequality
$$\sup_{m\in\N}\int_{B(0, 10)}|D(u_m\circ l_m)(x)|^2\dx<\infty.$$
It implies that $\{h_m\}$ is a bounded and (possibly after taking a subsequence) weakly convergent sequence in the space $W^{1, 2}(B(0, 10),\ B(0, 10))$, $h_m\wto h$, where $h$ is the pointwise limit of $h_m$. 

\subsection{\texorpdfstring{ Degree satisfies $\mathbf{\Deg(h,B(0, r), y )=2}$}{Deg(h, B(0,r), y)=2} } \ \\
We claim that for every point $y\in B\big((\tfrac{1}{2}, 0, 0), \tfrac{1}{2}\big)$ and for every radius ${r\in(\tfrac{5}{16}, \tfrac{7}{16})}$ we have 
$$ \Deg\left(h, B(0, r), y\right)=2. $$
Indeed, for every $r\in(\tfrac{5}{16}, \tfrac{7}{16})$ it holds that $h\in W^{1, 2}(\partial B(0, r))\cap L^\infty(\partial B(0, r))$. Fix such an $r$. 
The mappings $h_m$ map the sphere $\partial B(0, r)$ onto three bubbles, see Fig.~\ref{fig:overview}. In the limit, the filled and the dashed bubbles become topological spheres with the same orientation as the original sphere $B(0,r)$, and the dotted bubble disappears, see Fig.~\ref{fig:limit}.
Therefore it is not difficult to see that the degree of $h$ is $2$ inside the smaller topological sphere.

\section{For other degrees}\ \\
In this section, we explain an idea of how to construct a sequence of bounded and weakly convergent homeomorphisms $\{h_m\}\subset W^{1, 2}(B(0, 10), B(0, 10))$ which are identity on the boundary $\partial B(0, 10)$, such that the weak limit has degree $k\in \Z$ on a subset of positive measure. For $k=0, -1, 1$, the original construction by Conti and De Lellis in~\cite{CD2003} already gives the desired result. For $k=2$ please see our construction above.  

\begin{figure}[H]
	\begin{center}\input{image9.tex} \end{center}
	\caption{The construction of arbitrary degree. On the left-hand side for degree $3$, on the right-hand side for degree -2. The dotted lines will disappear in the limit, the dashed will remain.}
	\label{fig:higherDegree}
\end{figure}
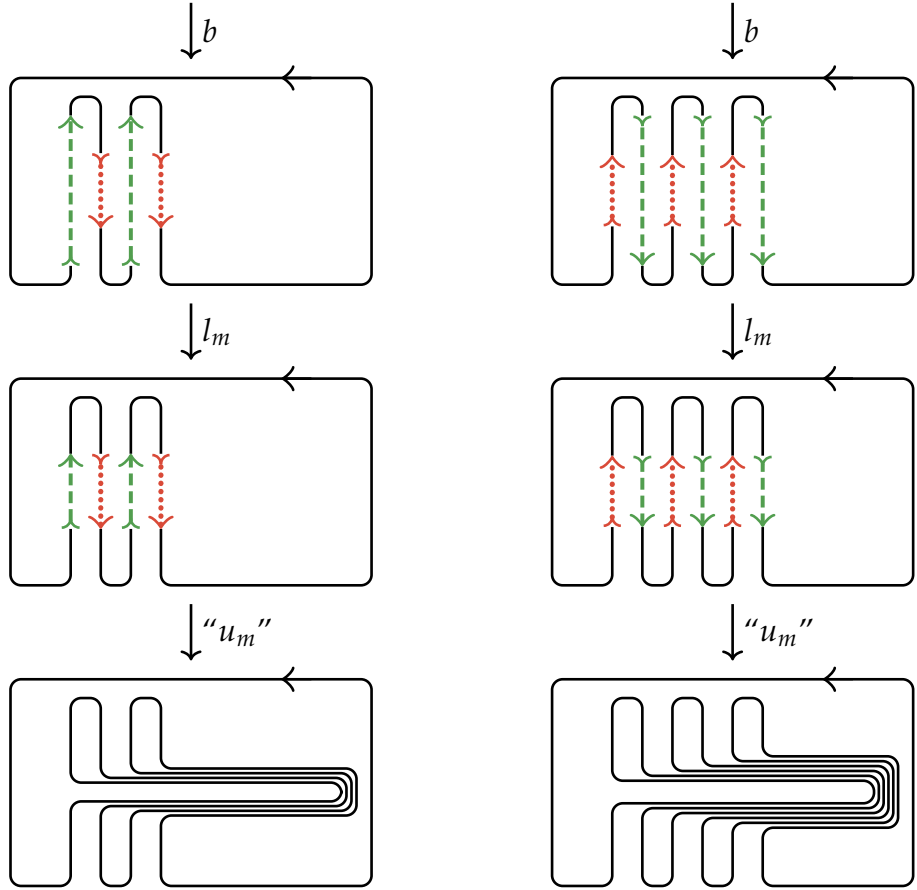

For degrees $|k|\geq 2$ the construction is similar to the case $k=2$. We need to modify our mapping $b$ and $l_m$ from Fig.\ \ref{fig:wigglycircle}. Instead of three, we need to create the appropriate number of bubbles to achieve the desired degree $k$. To do that we use a bi-Lipschitz mapping $b$ that maps the sphere onto a sufficiently wiggly shape as in Fig.\ \ref{fig:higherDegree}. And by choosing which arrows to shrink using $l_m$ we can decide which bubbles will disappear in the limit. Arrows that are shrunk before applying the \qm{bubble making function $u_m$} will not disappear and can change the final degree, the arrows that do not change size by $l_m$ will disappear in the limit and so will not affect the final degree.

\bibliographystyle{plain}
\bibliography{Bibliography} 

\end{document}

%% file: HeadTikz.tex
\usepackage{tikz}
\usetikzlibrary{patterns}
\usetikzlibrary {arrows.meta}

\newlength{\hatchspread}
\newlength{\hatchthickness}
\newlength{\hatchshift}
\newcommand{\hatchcolor}{}
\tikzset{hatchspread/.code={\setlength{\hatchspread}{#1}},
         hatchthickness/.code={\setlength{\hatchthickness}{#1}},
         hatchshift/.code={\setlength{\hatchshift}{#1}},
         hatchcolor/.code={\renewcommand{\hatchcolor}{#1}}}
\tikzset{hatchspread=3pt,
         hatchthickness=0.4pt,
         hatchshift=0pt,
         hatchcolor=black}
\pgfdeclarepatternformonly[\hatchspread,\hatchthickness,\hatchshift,\hatchcolor]
   {vypln}
   {\pgfqpoint{\dimexpr-2\hatchthickness}{\dimexpr-2\hatchthickness}}
   {\pgfqpoint{\dimexpr\hatchspread+2\hatchthickness}{\dimexpr\hatchspread+2\hatchthickness}}
   {\pgfqpoint{\dimexpr\hatchspread}{\dimexpr\hatchspread}}
   {
    \pgfsetlinewidth{\hatchthickness}
    \pgfpathmoveto{\pgfqpoint{0pt}{\dimexpr\hatchspread+\hatchshift}}
    \pgfpathlineto{\pgfqpoint{\dimexpr\hatchspread+0.15pt+\hatchshift}{-0.15pt}}
    \ifdim \hatchshift > 0pt
      \pgfpathmoveto{\pgfqpoint{0pt}{\hatchshift}}
      \pgfpathlineto{\pgfqpoint{\dimexpr0.15pt+\hatchshift}{-0.15pt}}
    \fi
    \pgfsetstrokecolor{\hatchcolor}
    \pgfusepath{stroke}
   }

   \pgfdeclarepatternformonly[\hatchspread,\hatchthickness,\hatchshift,\hatchcolor]
   {vyplnRotated}
   {\pgfqpoint{\dimexpr-2\hatchthickness}{\dimexpr-2\hatchthickness}}
   {\pgfqpoint{\dimexpr\hatchspread+2\hatchthickness}{\dimexpr\hatchspread+2\hatchthickness}}
   {\pgfqpoint{\dimexpr\hatchspread}{\dimexpr\hatchspread}}
   {
    \pgfsetlinewidth{\hatchthickness}
    \pgfpathmoveto{\pgfqpoint{0pt}{\hatchshift}}
    \pgfpathlineto{\pgfqpoint
		{\dimexpr\hatchspread+0.15pt}{\dimexpr\hatchspread+0.15pt+\hatchshift}
	}
    \ifdim \hatchshift > 0pt
      \pgfpathmoveto{\pgfqpoint{0pt}{\hatchshift}}
      \pgfpathlineto{\pgfqpoint{\dimexpr0.15pt+\hatchshift}{-0.15pt}}
    \fi
    \pgfsetstrokecolor{\hatchcolor}
    \pgfusepath{stroke}
   }
      \newcommand{\delkaCarky}{0.1}
      \newcommand{\carkaX}[1]{\draw (#1,-\delkaCarky)--(#1,\delkaCarky);}
      \newcommand{\carkaY}[1]{\draw (-\delkaCarky,#1)--(\delkaCarky,#1);}

      \definecolor{myGreen}{rgb}{0.32,0.62,0.31}
      \definecolor{myRedIsh}{rgb}{0.673333333,0.4,0.25}
      \definecolor{myGreenIsh}{rgb}{0.496666667,0.51,0.28}
      \definecolor{myRed}{rgb}{0.85,0.29,0.22}
	  \definecolor{myGreenRed}{rgb}{0.585,0.455,0.265}

      \definecolor{myGreenBlue}{rgb}{0.375,0.615,0.53}

      \definecolor{myOrange}{rgb}{0.8,0.47,0.2}
      \definecolor{myYellow}{rgb}{1.0,0.78,0.43}
      \definecolor{myLightGreen}{rgb}{0.65,0.76,0.38}
      \definecolor{myBlue}{rgb}{0.43,0.61,0.75}

	  \newcommand{\openArea}[3][3pt]{
		\path [draw=#2!50,fill=#2!10] #3;
		\path [draw=none,fill=myRed,
	semitransparent,pattern=vypln,hatchspread=8pt,
	hatchthickness=0.5pt,hatchcolor=#2,hatchspread=#1] #3;	
	\path [draw=#2,thick,dashed] #3;
		}

		\newcommand{\openAreaR}[3][3pt]{
			\path [draw=#2!50,fill=#2!10] #3;
			\path [draw=none,fill=myRed,
		semitransparent,pattern=vyplnRotated,hatchspread=8pt,
		hatchthickness=0.5pt,hatchcolor=#2,hatchspread=#1] #3;	
		\path [draw=#2,thick,dashed] #3;
			}
	
	\newcommand{\labelOverBackground}[3]{
		\node[text=#1,inner sep = 2pt,rounded corners=3pt,
	draw, fill, color=#1!10] at #2 {#3};
	\node[text=#1] at #2 {#3};
	}
	\usetikzlibrary{shadows.blur}
	\makeatletter
	\colorlet{pstb@shadow@color}{red!\pgfbs@opacity!red}%
	\makeatother

	\newcommand{\labelX}[3][black]{
		\node[anchor=north, fill, color=white,opacity=0.8, text opacity=1, text=#1, rounded corners=3pt, inner sep=1.5, outer sep = 2.5] at (#2,0) {$#3$};
		\carkaX{#2}
	}
	\newcommand{\labelY}[3][black]{
		\node[anchor=east, fill, color=white,opacity=0.8, text opacity=1, text=#1, rounded corners=3pt, inner sep=1.5, outer sep = 2.5] at (0,#2) {$#3$};
		\carkaY{#2}
	}

%% file: image4.tex
\begin{tikzpicture}[scale=1.7]
 	\clip  (-1.5,1.7) rectangle (6.7,-4);
	\newcommand{\sirka}{1.75}
	\newcommand{\delkasipky}{0.09}
	\newcommand{\arrowLeftImageFour}[3]{\draw[line width=\sirka/1.3333333,#3]
	(#1+\delkasipky,#2+\delkasipky)--(#1,#2)--
	(#1+\delkasipky,#2-\delkasipky);}
	\newcommand{\arrowRightImageFour}[3]{\draw[line width=\sirka/1.333333,#3]
	(#1-\delkasipky,#2+\delkasipky)--(#1,#2)--
	(#1-\delkasipky,#2-\delkasipky);}
	\newcommand{\dimage}{1.75}
	\begin{scope}
		\newcommand{\uhelimg}{20}
		\draw[dashed, myRed, line width=\sirka] (-\uhelimg:1) arc(-\uhelimg:\uhelimg:1);
		\draw[line width=\sirka, myGreen] (\uhelimg:1) arc(\uhelimg:360-\uhelimg:1);
		\node[anchor=south west] at (45:0.9) {$\partial B(0,r)$};
		\node[anchor=east, myRed] at (0:1) {$C_{m,r}$};
		\arrowRightImageFour{0}{1}{myGreen}
		\arrowLeftImageFour{0}{-1}{myGreen}
	\end{scope}
	\draw[line width = 1.5pt, ->] (1.4,0) -- (1+\dimage-0.4,0);
	\node[anchor=south] at (1+\dimage/2,0) 
		{$u_m$};	
	\begin{scope}[shift={(2+\dimage,0)},yscale=0.7]
		\path[smooth, draw=myGreen, line width=\sirka] 
		(-1,-0.1)
		.. controls (2.5,-2.5) and (2.5,2.5) ..
		(-1,0.1);
		\path[dashed, smooth, draw=myRed, line width=\sirka] 
		(-1,-0.1)
		.. controls (3.5,-4) and (3.5,4) ..
		(-1,0.1);
		\node at (0.35,0) {$\deg(u_m, B(0,r),y)=0$};
		\node[myGreenRed] at (1.9,1.5) {$u_m(\partial B(0,r))$};
		\arrowLeftImageFour{0.7}{0.77}{myGreen}
		\arrowRightImageFour{0.7}{-0.77}{myGreen}
		\arrowLeftImageFour{1.2}{-1.2}{myRed}
		\arrowRightImageFour{1.2}{1.2}{myRed}
	\end{scope}
	\draw[line width = 1.5pt, ->] (2.9,-0.7) .. controls (2.6,-1) and (2.4,-1.2) .. (2.4,-1.6);
	\node[anchor=south] at (2.9636,-1.3424) {$m\to\infty$};
	\begin{scope}[shift={(1+\dimage/2-0.5,-2.7)},yscale=0.7]
		\path[smooth, draw, line width=\sirka]
		(-1,0)
		.. controls (3,-3.25) and (3,3.25) ..
		(-1,0) -- cycle;
		\node at (0.6,0) {$\deg(u, B(0,r),y)=-1$};
		\node at (2.15,1) {$u(\partial B(0,r))$};
		\arrowLeftImageFour{0.97}{0.94}{black}
		\arrowRightImageFour{0.97}{-0.94}{black}
	\end{scope}
\end{tikzpicture}

%% file: image5.tex
\begin{tikzpicture}[scale=2.3]
\newcommand{\delkasipky}{0.05}
\newcommand{\arrowLeft}[3]{\draw[line width=0.75,#3]
(#1+\delkasipky,#2+\delkasipky)--(#1,#2)--
(#1+\delkasipky,#2-\delkasipky);}
\newcommand{\arrowRightImageFive}[3]{\draw[line width=0.75,#3]
(#1-\delkasipky,#2+\delkasipky)--(#1,#2)--
(#1-\delkasipky,#2-\delkasipky);}
\draw[myRed, line width=1]
plot [smooth cycle, tension=0.6] coordinates 
{(0,0) (0.5,0.2) (1.7,0.2) (2,0) (1.7,-0.2) (0.5,-0.2) (0,0)
(-0.1,0.425) (0.3,0.7) (1.9,0.7) (2.28,0) (1.9,-0.7) (0.3,-0.7) (-0.1,-0.425)};
\arrowRightImageFive{1}{0.225}{myRed}
\arrowRightImageFive{1}{0.757}{myRed}
\arrowLeft{1}{-0.225}{myRed}
\arrowLeft{1}{-0.757}{myRed}
\node at (1,0) {$\deg(h,B(0,r),y)=2$};
\node at (1.1,0.45) {$\deg(h,B(0,r),y)=1$};
\node at (1.2,0.9) {$\deg(h,B(0,r),y)=0$};
\node[myRed] at (2.5868,-0.6302) {$h(\partial B(0,r))$};
\node[opacity=0] at (-0.5868,-0.6302) {$h(\partial B(0,r))$};
\end{tikzpicture}

%% file: image7.tex
\begin{tikzpicture}[scale=1.7,
myDotted/.style={line cap=round, dash pattern=on 0pt off 1.5\pgflinewidth},
myDashed/.style={line cap=rect, dash pattern=on 1.5\pgflinewidth off 2\pgflinewidth}]
\newcommand{\distance}{0.7}
\clip  (-2-0.2,1.7121+0.2) rectangle (1+3+\distance+0.2,-1-\distance-2.7121-0.2);
\newcommand{\lineSirka}{1.3}
\newcommand{\arrowHook}{Hooks[arc=90,length=3.4,line width=1.275]}
\newcommand{\arrowHead}{Classical TikZ Rightarrow[length=3.8,line width=1.275]}
\begin{scope}
	\newcommand{\malyUhel}{40}
	\newcommand{\velkyUhel}{80}
	\draw[myGreen,myDashed,line width = \lineSirka * 1.1,arrows = {\arrowHook-\arrowHead}]
		(-\velkyUhel/2:1) arc(-\velkyUhel/2:\velkyUhel/2:1);
	\draw[black,line width = \lineSirka *  \lineSirka * 1] (\velkyUhel/2:1) arc(\velkyUhel/2:90-\malyUhel/2:1);
	\draw[myRed, myDotted, line width = \lineSirka *  \lineSirka * 1.1,arrows = {\arrowHook-\arrowHead}]
	(90-\malyUhel/2:1) arc(90-\malyUhel/2:90+\malyUhel/2:1);
	\draw[black,line width = \lineSirka * 1] (90+\malyUhel/2:1) arc(90+\malyUhel/2:360-\velkyUhel/2:1);
	\draw[black,line width = \lineSirka * 0.9,arrows = {-\arrowHead}] 
		(180:1) arc(180:225:1);
\end{scope}
\begin{scope}[shift={(2+\distance,0)}]
	\path[smooth, draw opacity=1, draw=black, line width = \lineSirka *  1]
(40:1)
arc(40:-180:1)
.. controls (-1,1.7321) and (0.8893,2.2168) .. (35:2)
arc(35:-20:2)
arc(-20:-20-180:0.45)
arc(-20:20:1.1)
.. controls (0.9008,0.7412) and (1.4754,0.7309) .. (40:1.6)
arc(40:220:0.3)
;
\newcommand{\malyUhel}{40}
\newcommand{\velkyUhel}{80}
\draw[white,line width = \lineSirka * 1.5] (-\velkyUhel/2:1) arc(-\velkyUhel/2:\velkyUhel/2:1);
\draw[white,line width = \lineSirka * 1.5] (-\malyUhel/2:1.1) arc(-\malyUhel/2:\malyUhel/2:1.1);
\draw[myGreen,myDashed, line width = \lineSirka * 1.1,arrows = {\arrowHook-\arrowHead}] (-\velkyUhel/2:1) arc(-\velkyUhel/2:\velkyUhel/2:1);
\draw[myRed,myDotted, line width = \lineSirka * 1.1,arrows = {\arrowHead-\arrowHook}] (-\malyUhel/2:1.1) arc(-\malyUhel/2:\malyUhel/2:1.1);
\draw[black,line width = \lineSirka * 0.9,arrows = {-\arrowHead}] 
		(180:1) arc(180:225:1);
\draw[black,line width = \lineSirka * 0.9,arrows = {-\arrowHead}] 
		(0:2) arc(0:22:2);
\end{scope}
\begin{scope}[shift={(-1,-2.7121-\distance)}]
	\path[smooth, draw opacity=1, draw=black, line width = \lineSirka *  1]
(40:1)
arc(40:-180:1)
.. controls (-1,1.7321) and (0.8893,2.2168) .. (35:2)
arc(35:-20:2)
arc(-20:-20-180:0.45)
arc(-20:20:1.1)
.. controls (0.9008,0.7412) and (1.4754,0.7309) .. (40:1.6)
arc(40:220:0.3)
;
\newcommand{\malyUhel}{40}
\newcommand{\velkyUhel}{40}
\draw[white,line width = \lineSirka * 1.5] (-\velkyUhel/2:1) arc(-\velkyUhel/2:\velkyUhel/2:1);
\draw[white,line width = \lineSirka * 1.5] (-\malyUhel/2:1.1) arc(-\malyUhel/2:\malyUhel/2:1.1);
\draw[myGreen,myDashed, line width = \lineSirka * 1.1,arrows = {\arrowHook-\arrowHead}] (-\velkyUhel/2:1) arc(-\velkyUhel/2:\velkyUhel/2:1);
\draw[myRed,myDotted, line width = \lineSirka * 1.1,arrows = {\arrowHead-\arrowHook}] (-\malyUhel/2:1.1) arc(-\malyUhel/2:\malyUhel/2:1.1);
\draw[black,line width = \lineSirka * 0.9,arrows = {-\arrowHead}] 
		(180:1) arc(180:225:1);
\draw[black,line width = \lineSirka * 0.9,arrows = {-\arrowHead}] 
		(0:2) arc(0:22:2);
\end{scope}
\begin{scope}[shift={(2+\distance,-2.7121-\distance)}]
	\path[smooth, draw opacity=1, draw=black, line width = \lineSirka *  1]
(40:1)
arc(40:-180:1)
.. controls (-1,1.7321) and (0.8893,2.2168) .. (35:2)
arc(35:-20:2)
arc(-20:-20-180:0.45)
arc(-20:20:1.1)
.. controls (0.9008,0.7412) and (1.4754,0.7309) .. (40:1.6)
arc(40:220:0.3)
;
\newcommand{\malyUhel}{40}
\newcommand{\velkyUhel}{40}
\draw[white,line width = \lineSirka * 1.5] (-\velkyUhel/2:1) arc(-\velkyUhel/2:\velkyUhel/2:1);
\draw[white,line width = \lineSirka * 1.5] (-\malyUhel/2:1.1) arc(-\malyUhel/2:\malyUhel/2:1.1);
\path[smooth, draw=myRed, myDotted, line width = \lineSirka * 1,arrows = {\arrowHook-\arrowHead}]
(20:1.1)
.. controls (1.2443,-0.2026) and (1.44,0.5) .. (1.7351,0.1528)
.. controls (1.8216,0.0509) and (1.8216,-0.0509) .. (1.7351,-0.1528)
.. controls (1.44,-0.5) and (1.2443,0.2026) .. (-20:1.1)
;
\path[smooth, draw=myGreen,myDashed, line width = \lineSirka * 1,arrows = {\arrowHead-\arrowHook}]
(20:1)
.. controls (1.1818,-0.3232) and (1.43,0.38) .. (1.6652,0.119)
.. controls (1.7367,0.0397) and (1.7367,-0.0397) .. (1.6652,-0.119)
.. controls (1.43,-0.38) and (1.1818,0.3232) .. (-20:1)
;
\draw[black,line width = \lineSirka * 0.9,arrows = {-\arrowHead}] 
		(180:1) arc(180:225:1);
\draw[black,line width = \lineSirka * 0.9,arrows = {-\arrowHead}] 
		(0:2) arc(0:22:2);
\end{scope}
\newcommand{\delkaSipky}{0.5}
\newcommand{\delkaVelkeSipky}{0.8}
\newcommand{\opravaSipka}{0.0255}
\draw[line width = \lineSirka *  1pt, ->] (1+0.5*\distance-0.5*\delkaSipky,0) -- ++(\delkaSipky,0);
\node[anchor=south] at (1+0.5*\distance-\opravaSipka,0) {$b$};
\draw[line width = \lineSirka *  1pt, ->] (1+\distance/2+\delkaVelkeSipky/2,-1-\distance/2+\delkaVelkeSipky/2) -- (1+\distance/2-\delkaVelkeSipky/2,-1-\distance/2-\delkaVelkeSipky/2);
\node[] at (1+\distance/2-0.15+\opravaSipka,-1-\distance/2+0.15+\opravaSipka) {$l_m$};
\draw[line width = \lineSirka *  1pt, ->] (1+0.5*\distance-0.5*\delkaSipky,-2.7121-\distance) -- ++(\delkaSipky,0);
\node[anchor=south] at (1+0.5*\distance-\opravaSipka,-2.7121-\distance) {$u_m$};
\end{tikzpicture}

%% file: image2.tex
\begin{tikzpicture}[scale=6.222]
\renewcommand{\delkaCarky}{0.01}
\newcommand{\barvaPop}{black!40}
\begin{scope}[shift={(0,0)}]
	\draw[<->]  (-1.5/8,0) -- (5/8,0);
	\draw[<->]  (0,-1.5/8) -- (0,5/8);
	\openAreaR{myRed}{(-1/8,2/8) -- (1/8,2/8) -- (1/8,4/8) -- (-1/8,4/8) -- cycle}
	\openArea{myGreen}{(2/8,1/8) -- (2/8,-1/8)--(4/8,-1/8)--(4/8,1/8)-- cycle}
	\node at (0,5/8) [anchor=east, \barvaPop] {$x_2$};
	\node at (5/8,0) [anchor=south, \barvaPop] {$x_1$};
	 \labelY[\barvaPop]{2/8}{\tfrac28}
	 \labelY[\barvaPop]{4/8}{\tfrac48}
	 \draw[dotted, thick,\barvaPop] (1/8,0) -- (1/8,2/8);
		\labelX[\barvaPop]{1/8}{\tfrac18}
		\labelX[\barvaPop]{2/8}{\tfrac28}
		\labelX[\barvaPop]{4/8}{\tfrac48}
	\carkaY{-1/8}
	\node at (0,-1/8) [text=\barvaPop,anchor = east] {$-\tfrac18$};
	\draw[dotted, thick, \barvaPop] (0,-1/8) -- (2/8,-1/8);
	\labelOverBackground{myRed}{(0,3/8) }{$C_1$}
	\labelOverBackground{myGreen}{(3/8,0)}{$K$}
	\node[circle, draw, fill=myRed,myRed, inner sep=0pt,
	minimum width=3pt] at (0,4/8) {};
	\node[circle, draw, fill=myGreen,myGreen, inner sep=0pt,
		minimum width=3pt] at (4/8,0) {};
\end{scope}
\draw[line width = 1.5pt, ->] (6/8,0) -- (7/8,0);
\node[anchor=south] at (51/70,0) 
	{$b$};
\begin{scope}[shift={(9/8,0)}]
	\draw[<->]  (-1/8,0) -- (8/8,0);
	\draw[<->]  (0,-2/8) -- (0,5/8);
	\openArea{myGreen}{(2/8,1/8) -- (2/8,-1/8)--(4/8,-1/8)--(4/8,1/8)-- cycle}
	\openAreaR{myRed}{(5/8,1/8) -- (5/8,-1/8) -- (7/8,-1/8) -- (7/8,1/8) -- cycle}
	\node at (0,5/8) [anchor=east, black!50] {$x_2$};
	\node at (8/8,0) [anchor=south, black!50] {$x_1$};
	\labelX[\barvaPop]{2/8}{\tfrac28}
	\labelX[\barvaPop]{4/8}{\tfrac48}
	\labelX[\barvaPop]{5/8}{\tfrac58}
	\labelX[\barvaPop]{7/8}{\tfrac78}
	 \draw[dotted, thick, \barvaPop] (0,1/8) -- (2/8,1/8);
	 \draw[dotted, thick, \barvaPop] (0,-1/8) -- (2/8,-1/8);
	 \labelY[\barvaPop]{1/8}{\tfrac18}
	 \labelY[\barvaPop]{-1/8}{-\tfrac18}
	\labelOverBackground{myGreen}{(3/8,0)}{$b(K)$}
	\labelOverBackground{myRed}{(6/8,0)}{$b(C_1)$}
	\node[circle, draw, fill=myRed,myRed, inner sep=0pt,
		minimum width=3pt] (bod) at (5/8,0) {};
	\node[circle, draw, fill=myGreen,myGreen, inner sep=0pt,
		minimum width=3pt] at (4/8,0) {};
\end{scope}
\end{tikzpicture}

%% file: image1.tex
\begin{tikzpicture}[scale=9]
\renewcommand{\delkaCarky}{0.01}
\newcommand{\barvaPop}{black!40}
\newcommand{\mmm}{2.6}
	\labelY[\barvaPop]{-1/8/\mmm}{-\tfrac1{8m}}
	\draw[<->]  (-1/8,0) -- (8/8,0);
	\draw[<->]  (0,-1.5/8) -- (0,1.7/8);
	\openArea{myBlue}{(1/8,1/8/\mmm) -- (2/8,1/8) -- (2/8,-1/8) -- (1/8,-1/8/\mmm) -- cycle}
	\openArea{myBlue}{(5/8,1/8/\mmm) -- (4/8,1/8) -- (4/8,-1/8) -- (5/8,-1/8/\mmm) -- cycle}
	\openAreaR{myRed}{(5/8,1/8/\mmm) -- (5/8,-1/8/\mmm)--(7/8,-1/8/\mmm)--(7/8,1/8/\mmm)-- cycle}
	\openArea{myGreen}{(2/8,1/8) -- (2/8,-1/8) -- (4/8,-1/8) -- (4/8,1/8) -- cycle}
	\node at (0,1.7/8) [anchor=east, black!50] {$x_2$};
	\node at (8/8,0) [anchor=south, black!50] {$x_1$};
	\labelX[\barvaPop]{2/8}{\tfrac28}
	\labelX[\barvaPop]{7/8}{\tfrac78}
	\labelY[\barvaPop]{1/8}{\tfrac18}
	 \draw[dotted, thick, \barvaPop] (0,-1/8/\mmm) -- (1/8,-1/8/\mmm);
	 \draw[dotted, thick, \barvaPop] (0,1/8) -- (2/8,1/8);
	\labelOverBackground{myBlue}{(1.5/8,0)}{$T^L_m$}
	\labelOverBackground{myGreen}{(3/8,0)}{$b(K)=K$}
	\labelOverBackground{myBlue}{(4.5/8,0)}{$T^R_m$}
	\labelOverBackground{myRed}{(6/8,0)}{$b(C_m)$}
	\labelX[\barvaPop]{4/8}{\tfrac48}
	\labelX[\barvaPop]{5/8}{\tfrac58}
	\labelX[\barvaPop]{7/8}{\tfrac78}
	\labelX[\barvaPop]{1/8}{\tfrac18}
\end{tikzpicture}

%% file: image6.tex
\begin{tikzpicture}[scale=1.5]
\begin{scope}
	\draw (0,1) arc(90:270:1);
	\draw (0,0.1) arc(90:270:0.1);
	\draw (0,1.5) -- (0,-1.5);
	\draw (0,0.1) -- (1,0.1);
	\draw (0,-0.1) -- (1,-0.1);
	\draw (1,1.5) -- (1,-1.5);
	\draw (1,1) arc(90:-90:1);
	\draw (1,0.1) arc(90:-90:0.1);
	\node at (-0.6,0) {\texttt{a}};
	\node at (1.6,0) {\texttt{e}};
	\node at (0.5,0) {\texttt{c}};
	\node at (0.5,0.1+1.4/2) {\texttt{d}};
	\node at (-0.8,1.1) {\texttt{b}};
	\node at (1.8,1.1) {\texttt{f}};
	\node at (-0.2,0.3) {\texttt{a'}};
	\draw (-0.05,0) -- (-0.2*0.6,0.3*0.6);
	\node at (1.2,-0.3) {\texttt{e'}};
	\draw (1.05,0) -- (1+0.2*0.8,-0.3*0.8);
\end{scope}
	\draw[line width = 1pt, ->] (2.3,0) -- (2.9,0);
	\node[anchor=south] at (5.2/2,0) {$u_m$};	
\begin{scope}[shift={(4,0)}]
\newcommand{\eps}{0.3}
	\draw (0,\eps) arc(90:-90:\eps);
	\draw (0,\eps) -- (-1,1+\eps);
	\draw (0,-\eps) -- (-1,-1-\eps);
	\draw (0,\eps) -- (0,1.5);
	\draw (0,-\eps) -- (0,-1.5);
	 \draw[] plot [smooth] coordinates 
		{(0., 0.6) (0.12585, 0.75417) (0.30025, 0.8746) (0.51209, 0.94627) (0.74578, 0.95818) (0.98272, 0.90466) (1.20315, 0.78606) (1.3881, 0.60888) (1.52138, 0.38526) (1.59113, 0.13184) (1.59113, -0.13184) (1.52138, -0.38526) (1.3881, -0.60888) (1.20315, -0.78606) (0.98272, -0.90466) (0.74578, -0.95818) (0.51209, -0.94627) (0.30025, -0.8746) (0.12585, -0.75417) (0., -0.6)};
 	\draw[] plot [smooth] coordinates 
		{(0., 0.6) (0.13398, 0.80287) (0.33188, 0.96673) (0.58005, 1.07184) (0.85896, 1.10359) (1.14511, 1.05415) (1.4134, 0.92342) (1.6397, 0.71924) (1.8033, 0.45666) (1.88909, 0.15653) (1.88909, -0.15653) (1.8033, -0.45666) (1.6397, -0.71924) (1.4134, -0.92342) (1.14511, -1.05415) (0.85896, -1.10359) (0.58005, -1.07184) (0.33188, -0.96673) (0.13398, -0.80287) (0., -0.6)};
	\draw[] plot [smooth] coordinates 
		{(0., 0.69) (0.14879, 0.89164) (0.3611, 1.05185) (0.62289, 1.151) (0.91424, 1.17462) (1.21133, 1.1151) (1.48875, 0.97265) (1.72212, 0.75539) (1.89054, 0.47875) (1.97878, 0.16397) (1.97878, -0.16397) (1.89054, -0.47875) (1.72212, -0.75539) (1.48875, -0.97265) (1.21133, -1.1151) (0.91424, -1.17462) (0.62289, -1.151) (0.3611, -1.05185) (0.14879, -0.89164) (0., -0.69)};
	\draw[] plot [smooth] coordinates 
		{(0., 0.69) (0.15692, 0.94035) (0.39273, 1.14398) (0.69085, 1.27657) (1.02742, 1.32003) (1.37371, 1.26459) (1.699, 1.11001) (1.97371, 0.86575) (2.17247, 0.55014) (2.27673, 0.18866) (2.27673, -0.18866) (2.17247, -0.55014) (1.97371, -0.86575) (1.699, -1.11001) (1.37371, -1.26459) (1.02742, -1.32003) (0.69085, -1.27657) (0.39273, -1.14398) (0.15692, -0.94035) (0., -0.69)};
	\draw[] plot [smooth] coordinates 
		{(0., 0.9) (0.10262, 1.10747) (0.24303, 1.3001) (0.41855, 1.47104) (0.62526, 1.61397) (0.85813, 1.72337) (1.11119, 1.79463) (1.37766, 1.82431) (1.65022, 1.8102) (1.92122, 1.75142) (2.18293, 1.64847) (2.42779, 1.50323) (2.64868, 1.31888) (2.83909, 1.09987) (2.99339, 0.85169) (3.10702, 0.5808) (3.17658, 0.29435) (3.2, 0.) (3.17658, -0.29435) (3.10702, -0.5808) (2.99339, -0.85169) (2.83909, -1.09987) (2.64868, -1.31888) (2.42779, -1.50323) (2.18293, -1.64847) (1.92122, -1.75142) (1.65022, -1.8102) (1.37766, -1.82431) (1.11119, -1.79463) (0.85813, -1.72337) (0.62526, -1.61397) (0.41855, -1.47104) (0.24303, -1.3001) (0.10262, -1.10747) (0., -0.9)};
	\node at (-0.6,0) {\texttt{b}};
	\node at (1.6/2,0) {\texttt{a}};
	\node at (3.5/2,0.021) {\texttt{a'}};
	\node at (2.14,0.021) {\texttt{e'}};
	\node at (2.745,0) {\texttt{e}};
	\node at (-0.267878, 0.946716) {\texttt{d}};
	\node at (3, 1.5) {\texttt{f}};
	\node at (1.945,-1.3) {\texttt{c}};
	\draw[thick,<-] (1.945,-0.05) -- (1.945,-1.15);
\end{scope}
\end{tikzpicture}

%% file: image3.tex
\begin{tikzpicture}[scale=6]
	\renewcommand{\delkaCarky}{0.01}
	\newcommand{\barvaPop}{black!40}
	\begin{scope}[shift={(1,0)}, scale=0.5] 
		\path[fill=black!15] (0,1) arc(90:270:1);
		\path[fill=black!15] (1,1) arc(90:-90:1);
	\end{scope}
	\begin{scope}[shift={(0,0)},scale=0.444444] 
		\draw[<->]  (-4/8,0) -- (5/8,0);
		\draw[<->]  (0,-4/8) -- (0,5/8);
		\openAreaR{myRed}{(-1/8,2/8) -- (1/8,2/8) -- (1/8,4/8) -- (-1/8,4/8) -- cycle}
		\openArea{myGreen}{(2/8,1/8) -- (2/8,-1/8)--(4/8,-1/8)--(4/8,1/8)-- cycle}
		\node at (5/8,0) [anchor=south, \barvaPop] {$x_1$};
		 \draw[dotted, thick,\barvaPop] (1/8,0) -- (1/8,2/8);
		 \labelX[\barvaPop]{1/8}{\tfrac18}
		 \labelX[\barvaPop]{2/8}{\tfrac28}
		 \labelX[\barvaPop]{4/8}{\tfrac48}
		 \labelOverBackground{myRed}{(0,3/8) }{$C_1$}
		 \labelOverBackground{myGreen}{(3/8,0)}{$K$}
		\labelY[\barvaPop]{2/8}{\tfrac28}
		 \node[circle, draw, fill=myRed,myRed, inner sep=0pt,
		 minimum width=2pt] at (0,4/8) {};
		\node[circle, draw, fill=myGreen,myGreen, inner sep=0pt,
			minimum width=2pt] at (4/8,0) {};
		\node at (0,5/8) [anchor=east, \barvaPop] {$x_2$};
		\labelX{-3/8}{-\tfrac38}
	\end{scope}
	\draw[line width = 1pt, ->] (0.1464,0.2625) .. controls (0.1605,0.3183) and (0.2812,0.3794) .. (0.3404,0.3604);
	\node[anchor=south] at (0.2265,0.3371) {$b$};
	\begin{scope}[shift={(8/8,0)},scale=0.5]
		\draw[<->]  (-1/8,0) -- (8/8,0);
		\draw[<->]  (0,-2/8) -- (0,5/8);
		\openArea{myGreen}{(2/8,1/8) -- (2/8,-1/8)--(4/8,-1/8)--(4/8,1/8)-- cycle}
		\openAreaR{myRed}{(5/8,1/8) -- (5/8,-1/8) -- (7/8,-1/8) -- (7/8,1/8) -- cycle}
		\node at (0,5/8) [anchor=east, black!50] {$x_2$};
		\node at (8/8,0) [anchor=south, black!50] {$x_1$};
		\labelX[\barvaPop]{2/8}{\tfrac28}
		\labelX[\barvaPop]{4/8}{\tfrac48}
		\labelX[\barvaPop]{5/8}{\tfrac58}
		\labelX[\barvaPop]{7/8}{\tfrac78}
		 \draw[dotted, thick, \barvaPop] (0,1/8) -- (2/8,1/8);
		 \draw[dotted, thick, \barvaPop] (0,-1/8) -- (2/8,-1/8);
		 \node[anchor=east, fill, color=white,opacity=0, text opacity=1, text=\barvaPop, rounded corners=3pt, inner sep=1.5, outer sep = 2.5]
			 at (0,1/8) {$\tfrac18$};
		\carkaY{1/8}
		\node[anchor=east, fill, color=white,opacity=0, text opacity=1, text=\barvaPop, rounded corners=3pt, inner sep=1.5, outer sep = 2.5]
			 at (0,-1/8) {$-\tfrac18$};
		\carkaY{-1/8}
		\node[circle, draw, fill=myRed,myRed, inner sep=0pt,
			minimum width=2pt] (bod) at (5/8,0) {};
		\node[circle, draw, fill=myGreen,myGreen, inner sep=0pt,
			minimum width=2pt] at (4/8,0) {};
	\end{scope}
	\newcommand{\mezikruzisirka}{0.125}
\begin{scope}[shift={(0,0)},scale=0.44444] 
	\path [draw=none,fill=myLightGreen, fill opacity = 0.8, even odd rule]
		(0,0) circle (3/8-\mezikruzisirka/2) (0,0) circle (3/8+\mezikruzisirka/2);
		\labelX{-3/8}{-\tfrac38}
\end{scope}
\begin{scope}[shift={(1,0)}, scale=0.5]
	\newcommand{\tempuhel}{30}
	\path[smooth, draw opacity=0.8, draw=myLightGreen, line width = 9.5]
	(30:3/8)
	arc(30:-30:3/8)
	.. controls (0.0248,-0.7071) and (0.6,-1.2) .. (0,-1.2) 
	arc(-90:-270:1.2)
	-- (1,1.2)
	arc(90:-90:1.2)
	.. controls (0.5,-1.2) and (1.1002,-0.7071) .. (0.8002,-0.1875) 
	arc(180+30:180-30:3/8)
	.. controls (0.9585,0.4617) and (0.7682,0.6) .. (0.5625,0.6) 
	.. controls (0.3568,0.6) and (0.1665,0.4617) .. (30:3/8) 
	;
	\node at (-0.5,0) {\texttt{a}};
	\node at (1.5,0) {\texttt{e}};
\end{scope}
	\end{tikzpicture}

%% file: image9.tex
\begin{tikzpicture}[xscale=4, yscale=2.5,
	myDotted/.style={line cap=round, dash pattern=on 0pt off 1.5\pgflinewidth},
	myDashed/.style={line cap=rect, dash pattern=on 2\pgflinewidth off 2.5\pgflinewidth}]
	\newcommand{\arrowHook}{Hooks[arc=90,length=4,line width=1]}
	\newcommand{\arrowHead}{Classical TikZ Rightarrow[length=5,line width=1]}
	\newcommand{\width}{1.2}
	\newcommand{\radius}{4}
	\begin{scope}
		\begin{scope}[shift={(0,+1.6)}]
		\draw[black, line width=1, rounded corners=\radius]
		(0,0)
		foreach \p in {0.2,0.4}
		{-- (\p,0) -- (\p,1)-- (\p+0.1,1) -- (\p+0.1,0)}
		-- (\width,0)
		-- (\width,1.1)
		-- (0,1.1)
		-- cycle
		;
		{
			\foreach \pp in {0.2,0.4}
			{
			\draw[white,line width = 1.8] (\pp,0.1)-- (\pp,0.9);
			\draw[myGreen,myDashed,line width = 1.5,arrows = {\arrowHook-\arrowHead}] 
				(\pp,0.1)-- (\pp,0.9);
			}	
		}
		\foreach \pp in {0.3,0.5}
		{
			\draw[white,line width = 1.8] (\pp,0.3)--(\pp,0.7);
			\draw[myRed,myDotted,line width = 2,arrows = {\arrowHead-\arrowHook}] 
				(\pp,0.3)-- (\pp,0.7);
		}	
		\draw[black,line width=1,arrows = {-\arrowHead}] (\width-0.2,1.1) -- (\width-0.3,1.1);
		\end{scope}
		\begin{scope}
			\draw[black, line width=1, rounded corners=\radius]
			(0,0)
			foreach \p in {0.2,0.4}
			{-- (\p,0) -- (\p,1)-- (\p+0.1,1) -- (\p+0.1,0)}
			-- (\width,0)
			-- (\width,1.1)
			-- (0,1.1)
			-- cycle
			;
			{
				\foreach \pp in {0.2,0.4}
				{
				\draw[white,line width = 1.8] (\pp,0.3)-- (\pp,0.7);
				\draw[myGreen,myDashed,line width = 1.5,arrows = {\arrowHook-\arrowHead}] 
					(\pp,0.3)-- (\pp,0.7);
				}	
			}
			\foreach \pp in {0.3,0.5}
			{
				\draw[white,line width = 1.8] (\pp,0.3)--(\pp,0.7);
				\draw[myRed,myDotted,line width = 2,arrows = {\arrowHead-\arrowHook}] 
					(\pp,0.3)-- (\pp,0.7);
			}	
		\draw[black,line width=1,arrows = {-\arrowHead}] (\width-0.2,1.1) -- (\width-0.3,1.1);
			\end{scope}
		\begin{scope}[shift={(0,-1.6)}]
			\newcommand{\bhor}{(\width-0.05)}
			\newcommand{\ahor}{(\bhor-0.05)}
			\newcommand{\Lhor}{0.2}
			\newcommand{\Rhor}{0.5}
		\newcommand{\ver}{0.45 - 0.35 * \p / 1.4+0.35 * 0.2 / 1.4}
		\newcommand{\hor}{\ahor+(\p-\Lhor)*(\bhor-\ahor)/(\Rhor-\Lhor)}
		\newcommand{\verpl}{0.45 - 0.35 * \p / 1.4+0.35 * 0.1 / 1.4}
		\newcommand{\horpl}{\ahor+(\p+0.1-\Lhor)*(\bhor-\ahor)/(\Rhor-\Lhor)}
		\draw[black, line width=1, rounded corners=\radius]
		(0,0)
		foreach \p in {0.2,0.4} %
		{-- (\p,0) -- (\p,\ver)-- ({\hor},\ver) -- ({\hor},{1-(\ver)}) -- (\p,{1-(\ver)}) -- (\p,1)
		-- (\p+0.1,1) 
		-- (\p+0.1,{1-(\verpl)}) -- ({\horpl},{1-(\verpl)}) -- ({\horpl},\verpl) -- (\p+0.1,\verpl) -- (\p+0.1,0)}
		-- (\width,0)
		-- (\width,1.1)
		-- (0,1.1)
		-- cycle
		;
		\draw[black,line width=0.8,arrows = {-\arrowHead}] (\width-0.2,1.1) -- (\width-0.3,1.1);
		\end{scope}
		\draw[line width = 1, ->] (\width/2,1.5+1.6) -- ++(0,-0.3);
		\node[anchor=west] at (\width/2,1.5-0.15+1.6) {$b$};
		\draw[line width = 1, ->] (\width/2,1.5) -- ++(0,-0.3);
		\node[anchor=west] at (\width/2,1.5-0.15) {$l_m$};
		\draw[line width = 1, ->] (\width/2,-0.1) -- ++(0,-0.3);
		\node[anchor=west] at (\width/2,-0.25) {\qm{$u_m$}};
	\end{scope}
	\begin{scope}[shift={(\width+0.6,0)}]
		\begin{scope}[shift={(0,+1.6)}]
		\draw[black, line width=1, rounded corners=\radius]
		(0,0)
		foreach \p in {0.2,0.4,0.6}
		{-- (\p,0) -- (\p,1)-- (\p+0.1,1) -- (\p+0.1,0)}
		-- (\width,0)
		-- (\width,1.1)
		-- (0,1.1)
		-- cycle
		;
		{
			\foreach \pp in {0.2,0.4,0.6}
			{
			\draw[white,line width = 1.8] (\pp,0.31)-- (\pp,0.69);
			\draw[myRed,myDotted,line width = 2,arrows = {\arrowHook-\arrowHead}] 
			(\pp,0.31)-- (\pp,0.69);
			}	
			}
			\foreach \pp in {0.3,0.5,0.7}
			{
				\draw[white,line width = 1.8] (\pp,0.1)-- (\pp,0.9);
				\draw[myGreen,myDashed,line width = 1.5,arrows = {\arrowHead-\arrowHook}] 
				(\pp,0.1)-- (\pp,0.9);
		}	
		\draw[black,line width=1,arrows = {-\arrowHead}] (\width-0.2,1.1) -- (\width-0.3,1.1);
		\end{scope}
		\begin{scope}
			\draw[black, line width=1, rounded corners=\radius]
			(0,0)
			foreach \p in {0.2,0.4,0.6}
			{-- (\p,0) -- (\p,1)-- (\p+0.1,1) -- (\p+0.1,0)}
			-- (\width,0)
			-- (\width,1.1)
			-- (0,1.1)
			-- cycle
			;
			{
			\foreach \pp in {0.2,0.4,0.6}
			{
			\draw[white,line width = 1.8] (\pp,0.31)-- (\pp,0.69);
			\draw[myRed,myDotted,line width = 2,arrows = {\arrowHook-\arrowHead}] 
			(\pp,0.31)-- (\pp,0.69);
			}	
			}
			\foreach \pp in {0.3,0.5,0.7}
			{
				\draw[white,line width = 1.8] (\pp,0.31)-- (\pp,0.69);
				\draw[myGreen,myDashed,line width = 1.5,arrows = {\arrowHead-\arrowHook}] 
				(\pp,0.31)-- (\pp,0.69);
		}	
		\draw[black,line width=1,arrows = {-\arrowHead}] (\width-0.2,1.1) -- (\width-0.3,1.1);

			\end{scope}
		\begin{scope}[shift={(0,-1.6)}]
			\newcommand{\bhor}{(\width-0.05)}
			\newcommand{\ahor}{(\bhor-0.08)}
			\newcommand{\Lhor}{0.2}
			\newcommand{\Rhor}{0.7}
		\newcommand{\aver}{0.44}
		\newcommand{\bver}{0.31}
		\newcommand{\Lver}{\Lhor}
		\newcommand{\Rver}{\Rhor}
		\newcommand{\hor}{\ahor+(\p-\Lhor)*(\bhor-\ahor)/(\Rhor-\Lhor)}
		\newcommand{\ver}{\aver+(\p-\Lver)*(\bver-\aver)/(\Rver-\Lver)}
		\newcommand{\horpl}{\ahor+(\p+0.1-\Lhor)*(\bhor-\ahor)/(\Rhor-\Lhor)}
		\newcommand{\verpl}{(\aver+(\p+0.1-\Lver)*(\bver-\aver)/(\Rver-\Lver))}
		\draw[black, line width=1, rounded corners=\radius]
		(0,0)
		foreach \p in {0.2,0.4,0.6} %
		{-- (\p,0) -- (\p,{\ver})-- ({\hor},{\ver}) -- ({\hor},{1-(\ver)}) -- (\p,{1-(\ver)}) -- (\p,1)
		-- (\p+0.1,1) 
		-- (\p+0.1,{1-(\verpl)}) -- ({\horpl},{1-(\verpl)}) -- ({\horpl},{\verpl}) -- (\p+0.1,{\verpl}) -- (\p+0.1,0)}
		-- (\width,0)
		-- (\width,1.1)
		-- (0,1.1)
		-- cycle
		;
		\draw[black,line width=0.8,arrows = {-\arrowHead}] (\width-0.2,1.1) -- (\width-0.3,1.1);
		\end{scope}
		\draw[line width = 1, ->] (\width/2,1.5+1.6) -- ++(0,-0.3);
		\node[anchor=west] at (\width/2,1.5-0.15+1.6) {$b$};
		\draw[line width = 1, ->] (\width/2,1.5) -- ++(0,-0.3);
		\node[anchor=west] at (\width/2,1.5-0.15) {$l_m$};
		\draw[line width = 1, ->] (\width/2,-0.1) -- ++(0,-0.3);
		\node[anchor=west] at (\width/2,-0.25) {\qm{$u_m$}};
	\end{scope}
\end{tikzpicture}